\documentclass[11pt]{amsart}

\makeatletter
\usepackage{amssymb}
\usepackage{latexsym}
\usepackage{amsbsy}
\usepackage{amsfonts}

\def\marginpar#1{\ignorespaces}

\newtheorem{theorem}[equation]{Theorem}
\newtheorem{proposition}[equation]{Proposition}
\newtheorem{lemma}[equation]{Lemma}

\theoremstyle{definition}
\newtheorem{remark}[equation]{Remark}

\numberwithin{equation}{section}

\def\AArm{\fam0 \rm}%
\newdimen\AAdi%
\newbox\AAbo%
\def\AAk#1#2{\setbox\AAbo=\hbox{#2}\AAdi=\wd\AAbo\kern#1\AAdi{}}%

\newcommand{\BBone}{{\ensuremath{{\AArm 1\AAk{-.8}{I}I}}}}

\def\eqref#1{(\ref{#1})}
\def\eqlabel#1{\def\@currentlabel{#1}}

\def\formula#1{\def\@tempa{#1}\let\@tempb\theequation\def\theequation{%
\hbox{#1}}\def\@currentlabel{(\theequation)}$$}
\def\endformula{\leqno\hbox{(\@tempa)}$$\@ignoretrue\let\theequation\@tempb}

\def\given{\hskip5\p@\relax\vrule\@width.4\p@\hskip5\p@\relax}

\newcommand{\open}[1]{%
\par\normalfont\topsep6\p@\@plus6\p@\trivlist\item[\hskip\labelsep\itshape#1%
\@addpunct{.}]\ignorespaces}

\DeclareRobustCommand{\close}[1]{%
  \ifmmode 
  \else \leavevmode\unskip\penalty9999 \hbox{}\nobreak\hfill
  \fi
  \quad\hbox{$#1$}}

\newlength{\toskip}

\def \R {{\mathbb R}}

\def \D {{\mathbb D}}
\def \L {{\mathbb L}}

\def \phi {\varphi}

\makeatother

\begin{document}
\date{\today}

\title[QUADRATIC TRANSPORTATION COST ...]{A CRITERION FOR TALAGRAND'S QUADRATIC TRANSPORTATION COST INEQUALITY.}

 \author[P. Cattiaux]{\textbf{\quad {Patrick} Cattiaux $^{\spadesuit}$ \, \, }}

\address{{\bf {Patrick} CATTIAUX},\\ Ecole Polytechnique, CMAP, F- 91128 Palaiseau cedex,
CNRS 756\\ and Universit\'e Paris X Nanterre, \'equipe MODAL'X, UFR SEGMI\\ 200 avenue de la
R\'epublique, F- 92001 Nanterre, Cedex.} \email{cattiaux@cmapx.polytechnique.fr}

 \author[A. Guillin]{\textbf{\quad {Arnaud} Guillin $^{\diamondsuit}$}}
\address{{\bf {Arnaud} GUILLIN},\\ CEREMADE \, Universit\'e  Paris IX Dauphine, F- 75775 Paris cedex, CNRS
7534.} \email{guillin@ceremade.dauphine.fr}

\maketitle
 \begin{center}
 \textsc{$^{\spadesuit}$ Ecole Polytechnique \quad and \quad Universit\'e Paris X}
\medskip

\textsc{$^{\diamondsuit}$ Universit\'e Paris IX}
 \end{center}

\begin{abstract}
We show that the quadratic transportation cost inequality $T_2$ is equivalent to both a Poincar\'e
inequality and a strong form of the Gaussian concentration property.  The main ingredient in the
proof is a new family of inequalities, called modified quadratic transportation cost inequalities in
the spirit of the modified logarithmic-Sobolev inequalities by Bobkov and Ledoux \cite{BL97}, that
are shown to hold as soon as a Poincar\'e inequality is satisfied.
\end{abstract}

\bigskip

\bigskip

\bigskip

\section{\bf Introduction, framework and main results.}\label{I}

Transportation inequalities recently deserved a lot of interest, especially in connection with the
concentration of measure phenomenon (see \cite{Led99}, \cite{Led02}). Links with others renowned
functional inequalities, in particular logarithmic-Sobolev inequalities, were also particularly
studied (see \cite{BG99}, \cite{OV00}, \cite{BGL}, \cite{Led02} ...), as no direct or tractable criteria were available for this kind of inequalities.
\medskip

Given a metric space $(E,d)$ equipped with its Borel $\sigma$ field, the $\L^p$ Wasserstein distance
between two probability measures $\mu$ and $\nu$ on $E$ is defined as
\begin{equation}\label{1.1}
W_p(\mu  , \nu) \, = \, \Big( \inf_{\pi} \, \int_{E \times E} \, d^p(x,y) \, \pi(dx,dy) \,
\Big)^{1/p} \, ,
\end{equation}
where $\pi$ describes the set of all coupling of $(\mu,\nu)$ , i.e. the set of all probability
measures on the product space with marginal distributions $\mu$ and $\nu$.
\medskip

A probability measure $\mu$ is said to satisfy the $T_p(C)$ transportation cost inequality if for
all probability measure $\nu$,
\begin{equation}\label{1.2}
W_p(\mu  ,  \nu) \, \leq \, \sqrt{2C \, H(\nu,\mu)} \, ,
\end{equation}
where $H(\nu,\mu)$ stands for the Kullback-Leibler information (or relative entropy), i.e.
\begin{eqnarray*}
H(\nu,\mu) \, = \, \int \,  \log \, (\frac{d\nu}{d\mu})  \, d\nu \quad \textrm{if } \, \nu\ll\mu
\quad ; \quad +\infty \, \textrm{ otherwise.}
\end{eqnarray*}

As shown by K. Marton (\cite{Mar96}), $T_1$ implies a Gaussian type concentration for $\mu$.

Let us briefly recall the general argument, we shall use later.

For any Borel set $A$ with measure $\mu(A) \geq 1/2$ introduce $A_r^c \, = \, \{x \, , \, d(x,A)
\geq r\}$ and $d\mu_A \, = \, \frac{\BBone_A}{\mu(A)} \, d\mu$. Set $B$ for $A_r^c$ and assume that
$W_1(\nu,\mu) \, \leq \, \varphi(H(\nu,\mu))$ for all $\nu$. Then
\begin{eqnarray}\label{1.3}
r \, \leq \, W_1(\mu_B,\mu_A) & \leq & W_1(\mu_B,\mu) \, + \, W_1(\mu,\mu_A) \\ & \leq &
\varphi(H(\mu_A,\mu)) \, + \, \varphi(H(\mu_B,\mu)) \nonumber \\ & = & \varphi\left(\log \,
\frac{1}{\mu(A)}\right) \, + \, \varphi\left(\log \, \frac{1}{\mu(A_r^c)}\right) \, .\nonumber
\end{eqnarray}
When $\varphi(u) \, = \, \sqrt{2C u}$ we immediately obtain $$\mu(A_r^c) \, \leq \, \exp \left( \, -
\, 1/2C \, \left(r \, - \, \sqrt{2C \, \log(\frac{1}{\mu(A)})}\right)^2 \, \right) \, .$$

Hence criteria for $T_1$ to hold are very useful. Such a criterion was first obtained by Bobkov and
G\"{o}tze (\cite{BG99} Theorem 3.1) and recently discussed by Djellout, Guillin and Wu (\cite{DGW}
Theorem 2.3) where the following is proved

\begin{theorem}\label{1.4} \, \cite{DGW} \quad
$\mu$ satisfies $T_1$ if and only if there exist $\varepsilon > 0$ and $x_0\in E$ such that
\begin{equation*}
(EI_\varepsilon(2)) \qquad \qquad \int_E \, e^{\varepsilon \, d^2(x,x_0)} \, \mu(dx) \, < \, +
\infty \, .
\end{equation*}
\end{theorem}
\medskip

Unfortunately $T_1$ is not well adapted to dimension free bounds, while $T_2$ is, as shown by
Talagrand (\cite{Tal96}). The first example of measure satisfying $T_2$ is the standard Gaussian
measure (\cite{Tal96}), for which $C=1$. When $E$ is a complete smooth Riemannian manifold of finite
dimension, with $d$ the geodesic distance and $dx$ the volume measure, Otto and Villani
(\cite{OV00}) have studied the $T_2$ property for absolutely continuous probability measures
(Boltzmann measures) $${\bf (B.M)} \quad \quad \mu(dx) \, = \, e^{-V(x)} \, dx \, ,$$ for $V \in
C^2(E)$ in connection with the logarithmic-Sobolev inequality. Their method was recently improved by
Wang (\cite{Wa03}) in order to skip the curvature assumption made in \cite{OV00}.
\smallskip

{\it In the sequel we shall assume that $\mu$ is a Boltzmann measure with $V \in C^3$, and that the
diffusion process built on $E$ with generator $L=1/2 \, div(\nabla) \, - \, 1/2 \, \nabla V.\nabla$
is non explosive.}

This is assumption (A) in \cite{Wa03}. Conditions for non explosion are known. Here are some among
others when $E=\R^d$:
\begin{itemize}
  \item $V(x) \rightarrow +\infty$ as $|x| \rightarrow +\infty$ and $|\nabla V|^2 - \Delta V$ is
  bounded from below,
  \item $x.\nabla V(x) \, \geq \, - a |x|^2 - b$ for some $a$ and $b$ in $\R$,
  \item $\int \, |\nabla V|^2 \, d\mu \, < \, +\infty$ .
\end{itemize}
For the first two see e.g. \cite{Ro99} p.26, for the third one see e.g. \cite{CL96}.
\smallskip

The first result is thus

\begin{theorem}\label{1.5} \, \cite{OV00}, \cite{BGL}, \cite{Wa03}, (also see \cite{Co02}) \quad
Let $\mu$ be as above with finite moment of order 2. If $\mu$ satisfies the logarithmic-Sobolev
inequality (L.S.I) $$\int \, g^2 \, \log (g^2) \, d\mu \, - \, \Big(\int \, g^2 \, d\mu\Big) \, \log
\Big(\int \, g^2 \, d\mu\Big) \, \leq \, 2C \, \int \, |\nabla g|^2 \, d\mu \, ,$$ for all smooth
$g$, then $\mu$ satisfies $T_2(C)$.
\end{theorem}
\medskip

A partial converse of Theorem \ref{1.5} is also shown in \cite{OV00} (Corollary 3.1), namely

\begin{theorem}\label{1.6} \, \cite{OV00}, \cite{BGL} \quad
Let $\mu$ be as above with finite moment of order 2, and $E=\R^n$. If $\mu$ satisfies $T_2(C)$ and
the curvature assumption $$\textrm{Hess}( V) \, \geq \, K \, Id_n$$ for some $K\in \R$, then $\mu$
satisfies a logarithmic-Sobolev inequality (with some new constant $\bar{C}$), provided $$1 \, + \,
K \, C \,
> \, 0 \, .$$
\end{theorem}
\medskip

The latter restriction is very important and has to be compared with Wang's results (\cite{Wa97b}
and \cite{Wa01}) telling that a logarithmic-Sobolev inequality holds provided the curvature
assumption above and the integrability condition $EI_{\varepsilon}(2)$ in Theorem \ref{1.4} hold
with $$\varepsilon \, + \, K \, > \, 0 \, .$$ In other words, according to Theorem \ref{1.4} and
Theorem \ref{1.6}, under the curvature assumption, log-Sobolev, $T_1(C_1)$, $T_2(C_2)$ are all
equivalent for appropriate constants $C_1$ and $C_2$. Whether this equivalence holds without
restrictions on the constants or not was left open by these authors.
\medskip

Let us recall that another
approach of Theorems \ref{1.5} and \ref{1.6} was introduced by Bobkov, Gentil and Ledoux
(\cite{BGL}). First of all the general Monge-Kantorovitch duality theory indicates that for $p\geq
1$,
\begin{equation}\label{1.7}
W_p^p(\nu,\mu) \, = \, \sup \, \Big(\int \, g \, d\nu \, - \, \int \, f \, d\mu\Big) \, ,
\end{equation}
where the supremum is running over all pairs $(f,g)$ of measurable and bounded functions satisfying
\begin{equation}\label{1.8}
g(x) \, \leq \, f(y) \, + \, d^p(x,y) \, ,
\end{equation}
for every $x \, , \, y \, \in E$. In the infimum-convolution notation of Maurey (\cite{Mau91}),
$$Qf(x) \, = \, \inf_{y \in E} \, \Big(f(y) \, + \, d^p(x,y)\Big)$$ achieves the optimal choice.
Defining $$Q_tf(x) \, = \, \inf_{y \in E} \, \Big(f(y) \, + \, \frac 1t \, d^2(x,y)\Big)$$ one thus
introduces a semi-group satisfying the Hamilton-Jacobi initial value problem. Relying some kind of
hypercontractivity of this semi-group to the logarithmic-Sobolev inequality, these authors obtain
both Theorems \ref{1.5} and \ref{1.6} (without any curvature assumption for \ref{1.5} improving Otto
and Villani result as and before Wang's result, also see \cite{OV01}). In particular, the following
originally due to Otto and Villani is elementary shown in  \cite{BGL} subsection 4.1

\begin{theorem}\label{1.9} \,  \quad
Let $\mu$ be as above. If $\mu$ satisfies $T_2(C)$ then $\mu$ satisfies the Poincar\'e (or spectral
gap) inequality (S.G.I) i.e. for all smooth $f$ , $$\qquad \textrm{Var}_{\mu}(f) \, \leq \, C \,
\int \, |\nabla f|^2 \, d\mu \, . $$
\end{theorem}
\medskip

This result gives us a first hint on what should be the difference between $T_1$ and $T_2$ as $T_1$ is well known to hold when $(S.G.I.)$ fails (see \cite{DGW}, Remark 2.4).

One aim of the present paper is to show that actually
\medskip

\begin{theorem}\label{1.10}\quad

Let $\mu$ be as above. Then $\mu$ satisfies $T_2$ if and only if $\mu$ satisfies some Poincar\'e
inequality, the integrability condition $EI_{\varepsilon}(2)$ of Theorem \ref{1.4} and the following
 property :

there exists some $a>e^{\frac 32}$ and some constant $c(a)$ such that for all $\nu = h \, \mu$ with
$H(\nu,\mu) \, \leq \, 1/2$  $${\mathbf{(Tronc)}} \qquad W_2^2(\nu_a,\nu) \, \leq \, c(a) \,
H(\nu,\mu) \, ,$$ where $\nu_a \, = \, (1/\nu(h\leq a)) \, h \, \BBone_{h\leq a} \, \mu $.
\end{theorem}
\medskip

An explicit upper bound of the constant of this $T_2$ inequality in terms of the constants arising in the Poincar\'e's inequality,  $EI_{\varepsilon}(2)$, choice of $a$ and $c(a)$ can be computed (and $c(a)$ being given optimized in $a$). We shall see that, furthermore, if $EI_{\varepsilon}(2)$ holds, (Tronc) is implied by the following
Variance-Entropy property $${\mathbf {(Var-Ent)}} \qquad \int \, d^2(x,x_0) \, \BBone_{h>a} \, d\nu \,
\leq \, D(a) \, H(\nu,\mu) \, ,$$ for $\nu$ as before.
\medskip

The proof of Theorem \ref{1.10} lies on the recent work by Wang \cite{Wa03}. The limitation to the
finite dimensional setting is due to the fact we want to use Otto-Villani coupling technique as in
section 2 of \cite{Wa03}. However, one expects that Theorem \ref{1.10} extends to infinite
dimensional settings, as  path spaces. Indeed Theorem \ref{1.5} is extended to this setting in
\cite{Wa03} section 5 by using finite dimensional approximation (also see the final section in
\cite{DGW} for an approach using Girsanov transform), and Monge-Amp\`ere theory was extended to this
setting by Feyel and Ustunel (\cite{FU02} and \cite{FU03}). This will not be studied here.
\medskip

The proof of Theorem \ref{1.10} splits into two parts. In section \ref{II} we shall show that
(S.G.I) implies some transportation inequality for measures $\nu$ with a bounded density. Actually
we prove an interpolation result between (S.G.I) and (L.S.I) through a family of inequalities
$I(\alpha)$ introduced by Latala and Oleszkiewicz (see \cite{LO00}) for $0\leq \alpha \leq 1$ ,
\begin{equation}\label{1.11}
\quad I(\alpha) \quad \quad \sup_{p\in [1,2)} \, \frac{\int \, f^2 \, d\mu \, - \, \left(\int \, f^p
\, d\mu\right)^{\frac 2p}}{(2 \, - \, p)^{\alpha}} \, \leq \, C(\alpha) \, \int \, |\nabla f|^2 \,
d\mu \, .
\end{equation}
It is easily seen that $I(0)$ is the Poincar\'e inequality and $I(1)$ reduces to the
logarithmic-Sobolev inequality. Our first result is the following

\begin{theorem}\label{1.12} \quad

Let $\mu$ be as above. If $I(\alpha)$ holds then for all $\nu$ such that $\parallel
\frac{d\nu}{d\mu}\parallel_{\infty} \, \leq \, K$ the following modified transportation inequality
holds $$W_2(\nu,\mu) \, \leq \, D(\alpha) \, \left(\log \, K\right)^{\frac{1-\alpha}{2}} \,
\sqrt{C(\alpha) \, H(\nu,\mu)} \, ,$$ where $$D(\alpha) \, = \, 16 \, \exp \left(\frac{1-\alpha}{2}
\, (1 \, - \, \log(1-\alpha))\right) \, .$$
\end{theorem}

Remark that the previous Theorem and Marton's trick allow to recover the concentration property
shown in \cite{LO00}. Indeed, recall \eqref{1.3} and remark that the interesting $K$ is given by $K
= 1/\mu(A_r^c)$. We immediately see that if $I(\alpha)$ holds, $\mu(A_r^c)$ behaves like $\exp \, (
- C \, r^{\frac{2}{2-\alpha}})$.

We refer to \cite{Wa03}, \cite{Wa03b}, \cite{BR03}, \cite{Ch02} and \cite{BCR} for more refined
results in connection with $I(\alpha)$.
\medskip

Another characterization of $I(0)$ (i.e. (S.G.I)) is obtained in \cite{BGL} section 5.3 in terms of
a mixed transportation cost $W_L$. It is almost immediate that for some constants $C$ and $D$ , $$C
\,  W_L \, \leq \, W_1 \, \leq \, D \, (W_L \, + \, W_L^{\frac 12}) \, .$$ It follows from Corollary
5.1 in \cite{BGL} that $$\textrm{(S.G.I)} \quad \Rightarrow \quad W_1(\nu,\mu) \, \leq D \,
(H(\nu,\mu) \, + \, H^{\frac 12}(\nu,\mu)) \, .$$ But the behavior of Wasserstein metrics for large
entropy is easily related to exponential integrability thanks to the following elementary lemma
proved in section \ref{III}

\begin{lemma}\label{1.13}
Assume that $\mu$ satisfies $EI_\varepsilon(p)$  for some $\varepsilon > 0$. There exists a constant
$C(\varepsilon)$ such that for all $\nu$ satisfying $H(\nu,\mu) \, \geq \, 1$ , $W_p^p(\nu,\mu) \,
\leq \, C(\varepsilon) \, H(\nu,\mu) \, .$

Here $EI_\varepsilon(p)$ is defined as in \ref{1.4} with $d^p$ instead of $d^2$.
\end{lemma}

The first consequence of Lemma \ref{1.13} combined with Theorem \ref{1.4}, is that the
transportation inequalities $T_2$ and $T_1$ are ``equivalent'' for large entropy. Since Marton's
method is essentially concerned with large entropy, $T_2$ cannot furnish a better concentration
result than $T_1$.

The second consequence is that $T_2$ is mainly (and surprisingly) concerned with small entropy. That
is why one can expect that the modified transportation inequality \ref{1.12} together with a small
entropy (so that the density cannot be too big except on a small set) will yield the statement in
Theorem \ref{1.10}. The proof will be given in section \ref{III}.
\medskip

At this point we shall mention that the proof of Lemma \ref{1.13} is using the trivial independent
coupling. We learned from F. Bolley and C. Villani \cite{BV03} that, using a less trivial coupling
in \cite{Vil}, this statement can be greatly improved, in particular

\begin{proposition}\label{1.14} \quad {\bf Bolley and Villani}
$$EI_\varepsilon(p) \quad \quad \Rightarrow \quad \quad W_p^p(\nu,\mu) \, \leq \, C(\varepsilon) \,
\big( H(\nu,\mu) \, + \, H^{\frac 12}(\nu,\mu) \big) \, \,   .$$
\end{proposition}

Since (S.G.I) implies $EI_\varepsilon(1)$ , this result for $p=1$ is stronger than the one we
already recalled. Bolley and Villani are then able to get back Theorem \ref{1.4}
i.e. $EI_\varepsilon(2)$ is equivalent to the transportation inequality $T_1$, but with some better constant than in \cite{DGW}.
\medskip

Section \ref{II} mainly contains the proof of Theorem \ref{1.12}. Section \ref{III} contains the
proofs of Lemma \ref{1.13}, Theorem \ref{1.10} and related topics. In particular, going back to the
proof of Theorem \ref{1.10}, one can see that the main term to be controlled is either
$W_2^2(\nu_a,\nu)$ (using (Tronc)) or the left hand side in (Var-Ent). Elementary computations allow
to control the later and show

\begin{theorem}\label{1.15} \quad Let $\mu$ be as above.
\begin{enumerate}
\item[(1)]  If $EI_{\varepsilon}(2)$ holds and $a>e^{\frac 32}$ there exists some constant $c(a)$ such
that for all $\nu$ with $H(\nu,\mu) \, \leq \, 1/2$  $$W_2^2(\nu,\mu) \, \leq \, W_2^2(\nu_a,\mu) \,
+ \, c(a) \, H(\nu,\mu) \, \log(1/H(\nu,\mu)) \, .$$
\item[(2)]  If $EI_{\varepsilon}(2)$ and Poincar\'e are satisfied, there exists some constant
$C$ such that $$W_2^2(\nu,\mu) \, \leq \, C \, \left(1 \, + \, \sqrt{\log^+(1/H(\nu,\mu))}\right) \,
H(\nu,\mu) \, .$$
\end{enumerate}
\end{theorem}
\medskip

Even if this last inequality is not dimension free, one may use the concavity of $x\to x\sqrt{\log^+ x}$ to get some tensorization over the dimension for $\mu^{\otimes n}$ which thus verifies the preceding inequality with constant $C(n)=C\sqrt{\log n}$ (see \cite{Mar96}, or \cite[Th.2.5]{DGW} for dependent sequences) to be compared to $C.n$ obtained with the sole $T_1$ inequality.

\bigskip

{\bf {Acknowledgments.}} \quad We wish to thank Fran{\c c}ois Bolley and Cedric Villani for
numerous and fruitful exchanges. Michel Ledoux and Liming Wu are also gratefully acknowledged for
their kind interest in this work.

\bigskip

\section{\bf Modified transportation inequalities.}\label{II}

\begin{proof} {\it  of Theorem \ref{1.12}.}

Let $\nu$ be a probability measure such that $h \, = \, \frac{d\nu}{d\mu}$ satisfies $0<\beta \leq
h(x) \leq K$. We assume first that $h \in \D$ i.e. is the sum of a constant and a $C^{\infty}$
function with compact support.

Let $P_t$ denotes the $\mu$ symmetric semigroup with generator $L=1/2 \, div(\nabla) \, - \, 1/2 \, \nabla V.\nabla$, and define $\mu_t \, = \, (P_t h) \mu$.

Our method relies on Otto-Villani's coupling \cite{OV00}, refined by Wang \cite{Wa03}, whose idea is the following: to provide a coupling between $\mu_t$ and $\mu_{t+s}$ as $\pi_s(dx,dy)=\mu_t(dx)\delta_{\phi_s(x)}(dy)$ where $\phi_s$ is the well defined unique (under our assumptions) solution of the p.d.e.
$$\frac{\rm d}{{\rm d}s}\phi_s=-\xi_{t+s}\circ\phi_s,\quad \phi_0=Id,S\ge0$$
with $\xi_{t+s}(x)=\nabla \log P_{t+s}h(x)$.

Then, according to Otto and Villani \cite{OV00}, Lemma 2 (or more exactly its proof), or Wang \cite{Wa03} section 3,
\begin{eqnarray}\label{2.1}
A \, = \, \frac{d^+}{dt} \, \left(- \, W_2(\mu_t,\mu)\right) & \leq & \limsup_{s\rightarrow 0^+} \,
\frac 1s \, W_2(\mu_t,\mu_{t+s}) \, \\ & \leq & 2 \, \left(\int \, |\nabla \sqrt{P_t h}|^2  \,
d\mu\right)^{\frac 12} \, . \nonumber
\end{eqnarray}
Using $I(\alpha)$ we obtain for all $1\leq p < 2$ ,
\begin{equation}\label{2.2}
A \, \leq \, \frac{2 \, \sqrt{C(\alpha) \, (2-p)^{\alpha}} \, \, \, \int  |\nabla \sqrt{P_t h}|^2 \,
d\mu}{\sqrt{1 \, - \, \left(\int \, (P_t h)^{\frac p2} \, d\mu\right)^{\frac 2p}}} \, .
\end{equation}

Now using a similar argument as in Lemma 3.1 in \cite{Wa03} or simply the fact that $\D$ is a nice
core for the diffusion semigroup, the following computation is rigorous

\begin{equation}\label{2.3}
\quad \frac{d}{dt} \, \left( 1 \, - \, \left(\int \, (P_t h)^{\frac p2} \, d\mu\right)^{\frac
1p}\right) \,  =  \, - \frac 12 \, \left(\int \, (P_t h)^{\frac p2} \, d\mu\right)^{\frac 1p - 1} \,
\, \int (P_t h)^{\frac p2 - 1} \, L P_t h \, d\mu \,
\end{equation}
\begin{eqnarray*}
 & = &  \frac 12 \, \left(\int \, (P_t h)^{\frac p2} \, d\mu\right)^{\frac 1p - 1} \, \, \int \, (\frac p2 - 1)
\, (P_t h)^{\frac p2 - 2} \, |\nabla P_t h|^2 \, d\mu \nonumber \\ & = & \frac 12 \, \left(\int \,
(P_t h)^{\frac p2} \, d\mu\right)^{\frac 1p - 1} \, \, \int \, (\frac p2 - 1) \, (P_t h)^{\frac p2 -
1} \, |\nabla \sqrt{P_t h}|^2 \, d\mu \nonumber \\ & \leq & 0 \, . \nonumber
\end{eqnarray*}

Here we have used $\int \, \left(\phi'(g) \, Lg \, + \, \phi''(g) \, |\nabla g|^2\right) \, d\mu \,
= \, 0$ , with $\phi(g)=g^{\frac p2 - 1}$ .
\smallskip

But since $h\leq K$, $P_t h \leq K$ hence according to \eqref{2.2} and \eqref{2.3}
\begin{equation}\label{2.4}
A \,  \leq \,  \frac{2 \, \sqrt{C(\alpha) \, (2-p)^{\alpha}} \, \, \, \int  |\nabla \sqrt{P_t h}|^2
\, \frac{K^{1 - \frac p2}}{(P_t h)^{1 - \frac p2}} \,  d\mu}{\sqrt{1 \, - \, \left(\int \, (P_t
h)^{\frac p2} \, d\mu\right)^{\frac 1p}} \, \sqrt{1 \, + \, \left(\int \, (P_t h)^{\frac p2} \,
d\mu\right)^{\frac 1p}}}
\end{equation}
\begin{eqnarray*}
& \leq & - \, \frac{4 \, \sqrt{C(\alpha) \, (2-p)^{\alpha}}}{\sqrt{1 \, - \, \left(\int \, (P_t
h)^{\frac p2} \, d\mu\right)^{\frac 1p}}} \, \, \frac{K^{1 - \frac p2}}{(1 - p/2) \, \big(\int (P_t
h)^{\frac p2} d\mu \big)^{\frac 1p - 1}} \, \frac{d}{dt} \, \left( 1 \, - \, \left(\int \, (P_t
h)^{\frac p2} \, d\mu\right)^{\frac 1p}\right) \nonumber \\ & \leq & 16 \, \sqrt{C(\alpha)} \,
(2-p)^{\frac{\alpha}{2} - 1} \, K^{1 - \frac p2} \, \left( - \frac{d}{dt} \, \sqrt{\left( 1 \, - \,
\left(\int \, (P_t h)^{\frac p2} \, d\mu\right)^{\frac 1p}\right)}\right) \, .
\end{eqnarray*}
For the latter inequality we have used $\int (P_t h)^{\frac p2} \, d\mu \, \leq \, 1$ .
\medskip

It remains to integrate in $t$. Since $I(\alpha)$ implies (S.G.I), we know that $P_t h$ goes to 1 in
$\L^2(\mu)$ as $t$ goes to infinity. Arguing as in \cite{Wa03} p.10, one can show that
$W_2(\mu_t,\mu)$ goes to 0 as $t$ goes to $\infty$, so that we have obtained

\begin{eqnarray}\label{2.5}
W_2(\nu,\mu) & \leq &  16 \, \sqrt{C(\alpha)} \, (2-p)^{\frac{\alpha}{2} - 1} \, K^{1 - \frac p2} \,
\sqrt{\left( 1 \, - \, \left(\int \, h^{\frac p2} \, d\mu\right)^{\frac 1p}\right)} \, \\ & \leq &
16 \, \sqrt{C(\alpha)} \, (2-p)^{\frac{\alpha}{2} - 1} \, K^{1 - \frac p2} \, \sqrt{\left( 1 \, - \,
\left(\int \, h^{\frac p2} \, d\mu\right)^{\frac 2p}\right)} \, . \nonumber
\end{eqnarray}

Now we shall use the two following elementary inequalities for $p\in [1,2)$:
\begin{itemize}
  \item $1 \, - \, u^{\frac 2p} \, \leq \, \frac 2p \, (1-u)$ for $u\in [0,1]$,
\medskip

  \item $\xi \, \log \xi \, + \, 1 \, - \, \xi \, \geq \, 0 $ for $\xi >0$ .
\end{itemize}

The latter yields $\log \xi^k \, \geq \, 1 \, - \, \xi^{-k}$, hence $\xi \, \log \xi^k \, \geq \,
\xi \, - \xi^{1-k}$ and finally for $k=1-\frac p2$, $(1-\frac p2) \, \xi \, \log \xi \, \geq \, \xi
\, - \, \xi^{\frac p2}$. We apply  this with $h(x)=\xi$, integrate with respect to $\mu$ and use the
former inequality in order to get
\begin{equation}\label{2.6}
1 \, - \, (\int h^{\frac p2} d\mu)^{\frac 2p} \, \leq \, \frac 2p \, (1-\frac p2) \, H(\nu,\mu) \, .
\end{equation}
Plugging \eqref{2.6} into \eqref{2.5} furnishes (using $p\geq 1$)

\begin{equation}\label{2.7}
W_2(\nu,\mu) \, \leq \, 16 \, \sqrt{C(\alpha)} \, (2-p)^{\frac{\alpha-1}{2}} \, K^{1 - \frac p2} \,
\sqrt{H(\nu,\mu)} \, .
\end{equation}

It is now enough to optimize in $p$. The optimal value is obtained for $2-p = \frac{1-\alpha}{\log
K}$, and a simple calculation yields the exact bound in Theorem \ref{1.12}.

It remains to extend the result to densities $h$ that are no more bounded away from 0, by using
standard tools.
\end{proof}
\medskip

One may ask whether this modified transportation inequality is dimension free. It does not seem so.
Actually the only kind of modified inequalities we are able to tensorize (following the induction
method in \cite{Tal96}) are the ones where we replace $(\log K)^{1-\frac{\alpha}{2}}$ by
$K^{\theta}$ for $\theta > 1/2$. For the concentration property, such a bound furnishes a polynomial
tail estimate for $\mu(A_r^c)$, precisely $(1/r)^{\frac{1}{\theta}}$ which is not really exciting.
\bigskip

\section{\bf Exponential integrability and the proof of Theorem \ref{1.10}.}\label{III}

We start this section by the proof of the elementary Lemma \ref{1.13} showing that the obstruction
for $T_2$ to hold is in a neighborhood of $\mu$. Notice that except for the conclusion (i.e Theorem
\ref{1.10} itself) all the intermediate results are available in a general metric space.

\begin{proof} {\it of Lemma \ref{1.13}.}

Introduce the Young function
\begin{equation}\label{3.1}
\tau(u) \, = \, u \, \log^+(u) \, ,
\end{equation}
and its Legendre conjugate function $\tau^*(v) \, = \, v \, \BBone_{v<1} \, + \, e^{v-1} \,
\BBone_{v\geq 1}$.

Among all possible coupling of $(\mu,\nu)$, the simplest one is the independent one i.e. if we
denote $h \, = \, \frac{d\nu}{d\mu}$ , $$\pi(dx,dy) \, = \, h(x) \, \mu(dx) \, \mu(dy) \, .$$
Accordingly
\begin{eqnarray*}
W_p^p(\nu,\mu) & \leq & \int \, d^p(x,y) \, h(x) \, \mu(dx) \, \mu(dy) \, \\ & \leq & 2 \, N_\tau(h)
\, N_{\tau^*}(d^p) \, ,
\end{eqnarray*}
where $N_\tau$ and $N_{\tau^*}$ are the gauge norms in the corresponding Orlicz spaces, the second
inequality being the classical H\"{o}lder-Orlicz inequality (see e.g. \cite{RR} for all concerned
with Orlicz spaces). Recall that the gauge norm for a general Young function $\psi$ is defined as
$$N_\psi(g) \, = \, \inf \, \{\lambda > 0 \, , \, \int \, \psi(g/\lambda)(x,y) \, \mu(dx) \, \mu(dy)
\, \leq \, 1\} \, ,$$ such that an easy convexity argument yields
\begin{equation}\label{3.2}
N_\psi(g) \, \leq \, \max \, \{1 \, , \, \int \, \psi(g) d\mu \otimes d\mu\} \, .
\end{equation}

In addition remark that $$\int \, h \, \log^+(h) \, = \, \int \, h \, \log(h) \, - \, \int_{h<1} \,
h \, \log(h) \, \leq \,  \int \, h \, \log(h) \, + \, 1/e \, .$$ Hence if $H(\nu,\mu) \, \geq \, 1$
, $$1 \, \leq \, \int \, h \, \log^+(h) \, \leq \, (1+1/e) \, H(\nu,\mu) \, ,$$ and according to
\eqref{3.2} and what precedes $$W_p^p(\nu,\mu) \, \leq \, 2(1+1/e) \, N_{\tau^*}(d^p) \, H(\nu,\mu)
\, .$$ Finally, thanks to $I_\varepsilon(p)$ , $ N_{\tau^*}(d^p) \, < \, +\infty$ and the result
follows.
\end{proof}
\medskip

One can improve the preceding result by showing that (up to the constant) it holds for $H(\nu,\mu)$
bounded away from 0. But as quoted in Proposition \ref{1.14} one can also get a precise bound for
the behavior of the Wasserstein distances when entropy goes to 0.
\medskip

\begin{proof} {\it of Theorem \ref{1.10}.}
We now proceed with the proof of Theorem \ref{1.10}. It breaks into several lemmata. According to
Lemma \ref{1.13} and \eqref{3.2} we may and will assume that $H(\nu,\mu)$ is small enough.

\begin{lemma}\label{3.3}
Let $\nu \, = \, h \, \mu$ be a probability measure. If $a>e$, then
\begin{enumerate}
\item[(1)] \quad  $H(\nu,\mu) \, \geq \, \big(1 \, - \, 1/\log a\big) \, \int_{h>a} \, h \, \log h \, d\mu$ ,
\item[]
\item[(2)] \quad  $\nu(h>a) \, \leq \, \big( 1 \, / \, (\log a \, - \, 1) \big) \, H(\nu,\mu)$ .
\end{enumerate}
\end{lemma}
\begin{proof}
Again we start with $u \, \log u \, + \, 1 \, - \, u \, \geq \, 0$ which yields $$\int_{h\leq a} \,
h \, \log h \, d\mu \, + \, 1 \, - \, \int_{h\leq a} \, h \, d\mu \, \geq \, 0 \, ,$$ hence
$$H(\nu,\mu) \, \geq \, \int_{h > a} \, h \, \log h \, d\mu \, - \, \nu(h > a) \, .$$ (2) follows
immediately since $ \log h \, > \, \log a$ on $\{h > a\}$ . For (1) we have $$\nu(h>a) \, \leq \,
\int_{h>a} \, \frac{\log h}{\log a} \, h \, d\mu \, = \, (1/\log a) \, \int_{h>a} \, h \, \log h \,
d\mu \, .$$
\end{proof}

Now we introduce a cut-off for $\nu$ i.e. if $a>0$ we define
\begin{equation}\label{3.4}
\nu_a \, = \, (1/\nu(h\leq a)) \, h \, \BBone_{h\leq a} \, \mu \, .
\end{equation}

\begin{lemma}\label{3.5}
Let $\nu \, = \, h \, \mu$ be a probability measure such that $H(\nu,\mu) \, \leq \, 1/2$. If
$a>e^{\frac 32}$ and $\nu_a$ is given by \eqref{3.4}, then $$H(\nu_a,\mu) \, \leq \, \left(1 +
\frac{1}{2(\log a - 3/2)} + \frac{2}{\log a - 1}\right) \, H(\nu,\mu) \, .$$
\end{lemma}
\begin{proof}
\begin{eqnarray*}
H(\nu_a,\mu) & = & \int \, \frac{h \, \BBone_{h\leq a}}{\nu(h\leq a)} \,
\log\left(\frac{h}{\nu(h\leq a)}\right) \, d\mu \\ & \leq & H(\nu,\mu) \, + \, \left((1/\nu(h\leq
a)) - 1\right) \, \int_{h\leq a} \, h \, \log h \, d\mu \\ & & - \, \log(\nu(h\leq a)) \, - \,
\int_{h>a} \, h \, \log h \, d\mu \\ & \leq & H(\nu,\mu) \, + \, \frac{\nu(h>a)}{\nu(h\leq a)} \,
H(\nu,\mu) \, - \, \log(1-\nu(h>a)) \, .\\
\end{eqnarray*}
But if $0\leq x \leq 1/2$ , $- \log(1-x) \, \leq \, 2 x$, hence according to \eqref{3.3}(2), if
$H(\nu,\mu) \, \leq \, 1/2$, $\log(1-\nu(h>a)) \, \leq \, (2/(\log a -1)) \, H(\nu,\mu)$ and
$$\frac{\nu(h>a)}{\nu(h\leq a)} \, \leq \, \frac{H(\nu,\mu)}{\log a - 1 - H(\nu,\mu)}$$ and we get
the desired result.
\end{proof}

We shall now proceed with the proof of an intermediate result : Poincar\'e, $EI_\varepsilon(2)$ and
(Var-Ent) imply $T_2$.
\smallskip

Recall the dual formulation of $W_2$ in \eqref{1.8} and \eqref{1.9} i.e. $$W_2^2(\nu_a,\mu) \, = \,
\sup \, \Big(\int \, g \, d\nu_a \, - \, \int \, f \, d\mu\Big) $$  $$ \textrm{ for $f$ and $g$ such
that for all $x$ and $y$} \quad g(x) \, \leq \, f(y) \, + \, d^2(x,y) \, .$$ Remark that in the
above formula we may add the same constant to both $f$ and $g$ so that we may assume that $\int  f
\, d\mu \, = \, 0$ . In this case, integrating with respect to $\mu(dy)$ the condition \eqref{1.9}
we have
\begin{eqnarray*}
g(x) & \leq & \int \, f \, d\mu \, + \, \int \, d^2(x,y) \, \mu(dy) \\ & \leq & 2 \, d^2(x,x_0) \, +
\, 2 \, \int \, d^2(y,x_0) \, \mu(dy)  \, = \, q_2(x) \, .
\end{eqnarray*}
 Now
\begin{eqnarray*}
\int \, g \, d\nu_a & = & \int \, g \, \frac{h \, \BBone_{h\leq a}}{\nu(h\leq a)} \, d\mu \\ & = &
\frac{1}{\nu(h\leq a)} \, \int \, g \, d\nu \, - \, \frac{1}{\nu(h\leq a)} \, \int \, g \, h \,
\BBone_{h>a} \, d\mu \\ & \geq & \frac{1}{\nu(h\leq a)} \, \int \, g \, d\nu \, - \,
\frac{1}{\nu(h\leq a)} \, \int \, q_2 \, h \, \BBone_{h>a} \, d\mu \, .
\end{eqnarray*}
Hence, since $\nu(h\leq a) \leq 1$ ,
\begin{equation}\label{3.6}
W_2^2(\nu,\mu) \, \leq \, W_2^2(\nu_a,\mu) \, + \, \int \, q_2 \, h \, \BBone_{h>a} \, d\mu \, .
\end{equation}
Recall that $q_2$ is the sum of a constant term and $2 \, d^2(x,x_0)$. So we have to control
\begin{equation}\label{3.7}
\int \, h \, \BBone_{h>a} \, d\mu \, = \, \nu(h>a) \, ,
\end{equation}
and
\begin{equation}\label{3.8}
\int \, d^2(x,x_0) \, h \, \BBone_{h>a} \, d\mu \, ,
\end{equation}
by some constant times $H(\nu,\mu)$. For \eqref{3.7} we may just use \eqref{3.3}(2), and for
\eqref{3.8} we may just use the hypothesis (Var-Ent) in Theorem \ref{1.10}. So applying successively
\eqref{3.6}, \eqref{3.7}, \eqref{3.8}, Theorem \ref{1.12} and Lemma \ref{3.5}, if $H(\nu,\mu) \,
\leq \, 1/2$
\begin{eqnarray*}
W_2^2(\nu,\mu) & \leq &  W_2^2(\nu_a,\mu) \, + \, c(a) \, H(\nu,\mu) \\ & \leq & D^2(0) \, C(0) \,
\log a \, H(\nu_a,\mu) \, + \, c(a) \, H(\nu,\mu) \\ & \leq & C(\alpha,a) \, H(\nu,\mu) \, .
\end{eqnarray*}

For $H(\nu,\mu) \geq 1$ we may use Lemma \ref{1.13}, and for $H(\nu,\mu) \in [1/2,1]$ we may use
\eqref{3.1} and  \eqref{3.2} and get $$W_2^2(\nu,\mu) \, \leq \, 2 \, N_{\tau^*}(d^2) \, \leq \, 4
\,  N_{\tau^*}(d^2) \, H(\nu,\mu) \, .$$

Hence we have proved that Poincar\'e, $EI_\varepsilon(2)$ and (Var-Ent) imply $T_2$, that is a
consequence of Theorem \ref{1.10}.

We did so because we shall use this method later to evaluate \eqref{3.8} when (Var-Ent) property
fails to hold. Furthermore, (Var-Ent) is well suited to study (Tronc).

Indeed, according to a well known result in mass transportation theory (see \cite{Vil} Proposition
7.10) if (Var-Ent) is satisfied, if $a>e^{3/2}$ and $H(\nu,\mu)\le 1/2$,
\begin{eqnarray*}
W_2^2(\nu,\nu_a) & \leq & C \, \int \, d^2(x,x_0) \, |1 \, - \, \frac{\BBone_{h\leq a}}{\nu(h\leq
a)}| \, d\nu \, ,\\ & \leq & C \, \left( \frac{\nu(h>a)}{\nu(h\leq a)} \, \int_{h\leq a} \,
d^2(x,x_0) \, d\nu \, + \, \int_{h> a} \, d^2(x,x_0) \, d\nu \right) \\ & \leq & C' \, H(\nu,\mu) \,
,
\end{eqnarray*}
according to Lemma \ref{3.3}, the H{\"o}lder-Orlicz inequality, $EI_{\varepsilon}(2)$ and (Var-Ent).
Hence (Tronc) is a consequence of (Var-Ent), provided $EI_{\varepsilon}(2)$ is satisfied.
\medskip

To finish, we proceed with the end of the proof of Theorem \ref{1.10}. For one way, it is enough to
write for $H(\nu,\mu) \leq \, 1/2$
\begin{eqnarray*}
W_2^2(\nu,\mu) & \leq &  2 \, W_2^2(\nu_a,\mu)  \, + \, 2 \, W_2^2(\nu_a,\nu) \\ & \leq & C' \,
H(\nu,\mu)
\end{eqnarray*}
according to the distance property of $W_2$, Theorem \ref{1.12}, Lemma \ref{3.5} and the (Tronc)
property.

Conversely we already know that $T_2$ implies both a Poincar\'e inequality and $EI_\varepsilon(2)$.
It remains to show that it also implies (Tronc). But if $H(\nu,\mu) \leq \, 1/2$ ,
\begin{eqnarray*}
W_2^2(\nu,\nu_a) & \leq & 2 \, W_2^2(\nu_a,\mu) \, + \, 2 \, W_2^2(\nu,\mu) \\ & \leq & 2 \, C \,
\left(H(\nu,\mu) + H(\nu_a,\mu)\right) \\ & \leq & C' \, H(\nu,\mu)
\end{eqnarray*}
according to the distance property, $T_2$ and Lemma \ref{3.5}.
\end{proof}
\medskip

To conclude this section we shall proceed with the proof of Theorem \ref{1.15}.

\begin{proof} {\it{of Theorem \ref{1.15}}}

{\it{Part (1).}} \quad According to \eqref{3.6}-\eqref{3.8} all we have to do is to estimate
$$\int_{h> a} \, d^2(x,x_0) \, d\nu \, .$$ Applying again the H{\"o}lder-Orlicz inequality and
$EI_\varepsilon(2)$ what we have to do is to estimate the Orlicz norm $$N_\tau(h \, \BBone_{h>a}) \,
,$$ i.e. we have to estimate $\lambda$ such that
\begin{equation}\label{3.9}
\int_{h>a} \, \frac{h}{\lambda} \, \log\left(\frac{h}{\lambda}\right) \, d\mu \, \leq \, 1 \, .
\end{equation}
According to Lemma \ref{3.3}, it is enough to have
\begin{equation}\label{3.10}
\frac{1}{\lambda} \, H(\nu,\mu) \, + \, \frac{1}{\lambda} \, \log\left(\frac{1}{\lambda}\right) \,
\nu(h>a) \, \leq \, 1 \, ,
\end{equation}
and it is easily seen that $\lambda \, \leq \, C(a) \, H(\nu,\mu) \, \log(1/H(\nu,\mu))$ .
\medskip

{\it{Part(2).}} \quad We shall be more accurate with the previous estimate. Indeed if
$EI_\varepsilon(2)$ and Poincar\'e are satisfied, it holds
\begin{eqnarray}\label{3.11}
W_2^2(\nu,\mu) & \leq & W_2^2(\nu_K,\mu) \, + \, \int_{h> K} \, d^2(x,x_0) \, d\nu \\ & \leq & C_1
\, \log\big(K/\nu(h\leq K)\big) \, H(\nu_K,\mu) \, + \, \int_{h> K} \, d^2(x,x_0) \, d\nu  \nonumber
\\ & \leq & C_2 \, \log(K) \, H(\nu,\mu) \, + \, \int_{h> K} \, d^2(x,x_0) \, d\nu \, , \nonumber
\end{eqnarray}
where we used successively \eqref{3.6}-\eqref{3.8}, Theorem \ref{1.12}, Lemma \ref{3.5} and previous
estimates (for small entropy, and large $K$).

Now we choose $K \, = \, 1/H^q(\mu,\nu)$ for some $q>0$ and we assume that $H(\nu,\mu)$ is small
enough (we already saw it is not a restriction). Lemma \ref{3.3}(2) furnishes
\begin{equation}\label{3.12}
\nu(h>K) \, \leq \,  H(\nu,\mu) / q \, \log(1/H(\nu,\mu)) \, ,
\end{equation}
so that the computation of $N_\tau(h \, \BBone_{h>K})$ as in \eqref{3.9}-\eqref{3.10} yields this
time $N_\tau(h \, \BBone_{h>K}) \, = \, C \, q^{-1} \, H(\nu,\mu)$ . Plugging this estimate into
\eqref{3.11} yields
\begin{equation}\label{3.13}
W_2^2(\nu,\mu) \, \leq \, \left(C_2 \, q \, \log(1/H(\nu,\mu)) \, + \, C_3 \, q^{-1}\right) \,
H(\nu,\mu) \, ,
\end{equation}
and the result follows optimizing in $q$ and using the same arguments as before for large entropy.
\end{proof}
\bigskip

\begin{remark}\label{3.14}

We hardly tried to improve the above estimates. For instance one can reduce the problem to estimate
$$\int_{1/H^q \, \geq \, e^{d^2} \, \geq h^p \, \geq  K^p} \, d^2(x,x_0) \, d\nu $$ for some fixed
$K$, $p>0$ large, $q>0$ small (this is left to the reader). Unfortunately we did not succeed in
removing the extra $\log(1/H(\nu,\mu))$ in this estimate, hence in Theorem \ref{1.15}. Actually we
do not know whether this is possible or not, only assuming Poincar\'e and the exponential
integrability. However we shall see in the next section that for some less general potentials $V$
one can do the job.
\end{remark}
\bigskip


\end{document}